\documentclass[twoside]{article}

\usepackage[preprint]{aistats2026}

\usepackage{amsmath,amssymb,amsthm,mathtools}
\usepackage{bbm,bm}
\usepackage{booktabs}
\usepackage{graphicx}
\usepackage{xcolor}
\usepackage{url}
\usepackage{microtype}
\usepackage{enumitem}
\usepackage{multirow}
\usepackage{algorithm, algpseudocode}
\usepackage{float}
\usepackage[round]{natbib}

\usepackage{hyperref}

\theoremstyle{plain}
\newtheorem{theorem}{Theorem}
\newtheorem{proposition}[theorem]{Proposition}
\newtheorem{lemma}[theorem]{Lemma}
\theoremstyle{definition}

\theoremstyle{remark}
\newtheorem{remark}[theorem]{Remark}
\newtheorem{corollary}[theorem]{Corollary}

\begin{document}

\twocolumn[
\aistatstitle{Matrix Phylogeny: Compact Spectral Fingerprints for Trap-Robust Preconditioner Selection}

\aistatsauthor{Jinwoo Baek}
\aistatsaddress{Department of Computer Science, Oregon State University}
]

\begin{abstract}
We introduce \emph{Matrix Phylogeny}, a framework for discovering family-level relationships among matrices via \textbf{Compact/Adaptive Spectral Fingerprints} (CSF/ASF). Our fingerprints are low-dimensional, eigendecomposition-free descriptors built from Chebyshev trace moments with Hutchinson sketches, following a simple affine rescaling to $[-1,1]$ that renders them invariant to permutation and (general) similarity transforms and robust to global scaling.

Across a comprehensive suite (E0--E6+), we observe \textit{phylogenetic compactness}: family identity is captured by only a handful of moments. A fixed small dimension already suffices—\textbf{CSF-$K=3\!-\!5$} achieves perfect clustering on four synthetic families and a five-family set including BA vs.~ER ($ARI = 1.0$, silhouettes~$\approx0.89$ at $K=3\!-\!5$) with the lowest runtimes, while \textbf{ASF} adapts on demand (\,$K^*$ mean~$8.4$, median~$9$\,) serving as a safety net. On a SuiteSparse mini-benchmark with Hutchinson ($p \approx 100$), both CSF-H and ASF-H reproduce the accuracy trend ($ARI = 1.0$). Against strong alternatives (eigenvalue histogram~+~Wasserstein, heat-kernel traces, WL-subtree), CSF-$K=5$ matches or exceeds accuracy while \emph{avoiding} eigendecompositions and using an order-of-magnitude smaller feature dimension (e.g., $K \le 10$ vs.~64/9153).

The descriptors are stable to noise (log--log slope $\approx1.03$, $R^2\!\approx\!0.993$) and enable a practical ``trap$\rightarrow$recommend'' pipeline for automated preconditioner selection. In an adversarial E6+ setting with a Probe-and-Switch mechanism, our physics-guided recommender attains near-oracle iteration counts (regret $\mathrm{p90}{=}0$), whereas a Frobenius 1-NN baseline exhibits large spikes ($\mathrm{p90}{\approx}34$--$60$). Altogether, CSF/ASF provide compact ($K{\le}10$), fast, and invariant fingerprints that unlock scalable, structure-aware search and recommendation over large matrix repositories; we recommend \textbf{CSF-$K{=}5$} by default and \textbf{ASF} when domain-specific adaptivity is desired.
\end{abstract}

\section{Introduction}
\label{sec:intro}

Large matrix repositories arise across simulation, optimization, and machine learning. Practitioners must group matrices by \emph{family resemblance} and make \emph{actionable} choices—e.g., selecting preconditioners—when a new matrix arrives. Entrywise distances ignore spectra; fully spectral methods are costly.

\textbf{Idea.} We propose \emph{Compact/Adaptive Spectral Fingerprints} (CSF/ASF): short descriptors that capture structure without computing eigenvectors or full spectra. After affine rescaling of the spectrum to $[-1,1]$, we take damped Chebyshev trace moments, yielding invariance to permutation, similarity transforms, and positive scaling. CSF fixes a tiny dimension ($K\!\in\!\{3,5\}$); ASF chooses $K$ via a simple stopping rule. Both support Hutchinson sketches. Pairwise distances induce a dendrogram—\emph{matrix phylogeny}—for retrieval, clustering, and recommendation across heterogeneous sizes.\footnote{By “eigendecomposition-free” we mean we never compute eigenvectors nor spectral histograms. For spectrum scaling we may compute scalar endpoints; on small dense cases our prototype obtains them via eigenvalue routines (no eigenvectors), while large runs can rely on matvec-only bounds (power/Lanczos or Gershgorin). Fingerprints are computed via matvec Chebyshev recurrences.}

\textbf{At a glance.} Across E0–E6 (synthetics, SuiteSparse, adversarial), a handful of moments identify families: CSF-$K{=}3$–$5$ attains ARI$=1.0$ on four synthetic families and on a five-family set including BA vs.\ ER (silhouette $\approx0.89$) at the lowest runtimes; ASF adapts conservatively ($K^*$ mean $8.4$, median $9$). On SuiteSparse with Hutchinson ($p\!\approx\!100$), CSF-H/ASF-H match this trend (ARI$=1.0$). Against strong alternatives (eigenvalue histograms + Wasserstein, heat-kernel traces, WL-subtree), CSF-$K{=}5$ matches or exceeds accuracy while using far fewer features ($K{\le}10$ vs.\ 64/9153) and avoiding full spectra. Fingerprints are noise-stable (slope $\approx\!1.03$, $R^2\!\approx\!0.993$). In adversarial E6+, a \emph{Probe-and-Switch} recommender guided by phylogenetic neighbors attains near-oracle iterations (regret $p_{90}\!=\!0$), whereas a Frobenius 1-NN baseline shows large spikes ($p_{90}{\approx}34$–$60$).

\textbf{Contributions.}
\begin{itemize}[leftmargin=2.0em]
  \item \textbf{Framework:} A practical pipeline for \emph{matrix phylogeny} via CSF/ASF that works across sizes without eigenvectors/full spectra.
  \item \textbf{Guarantees:} Invariance, Hutchinson concentration, noise stability, and adaptive stopping (formal statements in the appendix).
  \item \textbf{Evidence:} Strong accuracy at tiny dimensions, matching or exceeding spectral baselines that require full spectra.
  \item \textbf{Actionability:} A robust \emph{Probe-and-Switch} preconditioner selector with near-oracle iterations in adversarial settings.
\end{itemize}

\textbf{Organization.} Section~\ref{sec:method} details the method; Section~\ref{sec:theory} states guarantees; Section~\ref{sec:experiments} reports E0–E6; we conclude with limitations and related work.

\section{Preliminaries}
\label{sec:prelim}

\paragraph{Notation.}
For $A \in \mathbb{R}^{n \times n}$, let $\lambda_i(A)$ be eigenvalues, $\rho(A)=\max_i|\lambda_i(A)|$ the spectral radius, and $\mathrm{tr}(\cdot)$ the trace. For a polynomial $p$, $p(A)$ is the matrix polynomial.

\paragraph{Similarity-invariant spectral normalization.}
We compute fingerprints on $\widetilde{A}$ with spectrum in $[-1,1]$:
\begin{equation}
\label{eq:affine}
\begin{split}
\widetilde{A} &= \frac{A - mI}{r},\quad
m = \tfrac{1}{2}(\lambda_{\max} + \lambda_{\min}),\\
&r_0 = \tfrac{1}{2}(\lambda_{\max} - \lambda_{\min}),\\
&r = (1+\varepsilon_{\mathrm{rel}})\,r_0\quad(\varepsilon_{\mathrm{rel}}>0).
\end{split}
\end{equation}
For markedly non-Hermitian inputs we use radius scaling
$\widetilde{A}=A/\max(\rho(A),\varepsilon_{\mathrm{rad}})$ with tiny $\varepsilon_{\mathrm{rad}}>0$.
In practice, endpoints or $\rho(A)$ come from cheap bounds or a few power/Lanczos steps (Appendix~\ref{app:proofs}). Since the normalization depends only on spectral quantities and $T_k$ is a polynomial, $\mathrm{tr}\!\big(T_k(\widetilde A)\big)$ is invariant to permutations, similarities $S^{-1}AS$, and positive scalings $\alpha A$.

\paragraph{Eigendecomposition-free Chebyshev moments.}
With $T_0(x)=1$, $T_1(x)=x$, $T_{k+1}(x)=2xT_k(x)-T_{k-1}(x)$, define damped moments on $\widetilde{A}$:
\begin{equation}
\label{eq:moments}
s_k(A)=\mathrm{tr}\!\big(T_k(\widetilde{A})\big),\qquad
\widehat{s}_k(A)=e^{-\eta k}\, s_k(A).
\end{equation}
Hutchinson with $p$ probes $z_i$ ($\mathbb{E}[z_iz_i^\top]=I$) estimates traces:
\[
\widehat{\mathrm{tr}}(M)=\tfrac{1}{p}\sum_{i=1}^p z_i^\top M z_i,\quad
y_0=z_i,\ y_1=\widetilde{A}z_i,
\]
\[
y_{k+1}=2\widetilde{A}y_k-y_{k-1},
\]
so the cost is $O\!\big(pK\,\mathrm{cost}(\widetilde{A}\mathbf{v})\big)$.

\paragraph{Compact spectral fingerprints (CSF/CSF-H).}
For target $K$, form $d\in\mathbb{R}^K$ with
\[
d_0 := w_0\ \ (\text{default } w_0 = n \text{ in our experiments;}
\]
\[
\text{a size-neutral variant uses } w_0=1),
\]
and define the $\ell_2$-normalized fingerprint
\begin{equation}
\label{eq:csf}
\phi_K(A)=\frac{(d_0,\ldots,d_{K-1})}{\|(d_0,\ldots,d_{K-1})\|_2}\in\mathbb{R}^{K}.
\end{equation}
Exact traces give \emph{CSF}; Hutchinson gives \emph{CSF-H}. The adaptive variant (\emph{ASF}) chooses $K$ via a Hankel-ratio and an energy gate (Section~\ref{sec:method}).

\paragraph{Distances and evaluation.}
Unless noted, we use \textbf{Euclidean} distance between $\ell_2$-normalized fingerprints (cosine reported as an ablation; \emph{E6/E6+ uses cosine}). Clustering quality is measured by ARI and mean silhouette from the pairwise distance matrix.

\section{Methodology: Compact/Adaptive Spectral Fingerprints}\label{sec:method}

\paragraph{Goal.}
Given a (typically symmetric or symmetrizable) matrix $A\in\mathbb{R}^{n\times n}$, we seek a compact fingerprint $\phi(A)\in\mathbb{R}^{K_\star}$ that (i) is invariant to similarity/permutation/diagonal similarity transforms; (ii) is robust to noise and positive scaling; (iii) is computable in $\mathcal{O}(pK\,\mathrm{nnz}(A))$ using only sparse matvecs; and (iv) adapts $K_\star$ via data-driven stopping.

\subsection{Preliminaries and Notation}
Let $T_k$ be Chebyshev polynomials of the first kind ($T_0(x)=1$, $T_1(x)=x$, $T_{k+1}(x)=2xT_k(x)-T_{k-1}(x)$). We write $\mathrm{tr}(\cdot)$ for the trace, $\|\cdot\|_2$ for the spectral norm, and $\|\cdot\|_F$ for the Frobenius norm.

\paragraph{Affine normalization to $[-1,1]$.}
We use the affine map in Eq.~\eqref{eq:affine} with a small relative margin $\varepsilon_{\mathrm{rel}}$; for markedly non-Hermitian inputs we use radius scaling or symmetrization (see Section~\ref{sec:prelim}).

\paragraph{Implementation note.}
To obtain similarity-invariant scaling, we estimate only scalar spectral endpoints or a 2-norm bound—no eigenvectors or full spectra. On small dense inputs we may call an eigenvalue routine solely for endpoints (no eigenvectors); on large runs we use matvec-only bounds (power/Lanczos or Gershgorin). CSF/ASF themselves use matvec-based Chebyshev recurrences.

\subsection{Fixed-$K$ Chebyshev Spectral Fingerprints (CSF)}
For fixed $K\!\ge\!1$ and damping $\eta>0$, define
\[
\phi_K(A)=\bigl[d_0,\dots,d_{K-1}\bigr],\qquad d_k=e^{-\eta k}\,\mathrm{tr}\bigl(T_k(\widetilde A)\bigr).
\]
To avoid size leakage across different $n$, set
\[
d_0 \triangleq w_0,\quad \text{default } w_0=1\ \ (\text{size-neutral});
\]
the classical choice $w_0=n$ and very small $w_0\ll1$ are also supported.
\emph{We used $w_0=n$ in our experiments unless otherwise stated (see Section~\ref{sec:intro}).}
Finally, L2-normalize $\phi_K(A)\!\leftarrow\! \phi_K(A)/\|\phi_K(A)\|_2$.

\begin{algorithm}[H]
\caption{CSF-K (Fixed-$K$ Chebyshev Fingerprint; matrix form. At scale we use Hutchinson \emph{matvec} recurrences.)}
\label{alg:csf}
\begin{algorithmic}[1]
\Require $A$, target $K\!\ge\!1$, damping $\eta>0$, $w_0$ (default $=1$)
\State $\widetilde A \leftarrow \textsc{ScaleToUnitInterval}(A)$
\State $T_{\mathrm{prev}}\!\leftarrow I$, $T_{\mathrm{curr}}\!\leftarrow \widetilde A$
\State $d_0\!\leftarrow w_0$
\If{$K>1$} \State $d_1\!\leftarrow e^{-\eta}\, \mathrm{tr}(\widetilde A)$ \EndIf
\For{$k=2$ to $K-1$}
  \State $T_{\mathrm{next}}\!\leftarrow 2\widetilde A\,T_{\mathrm{curr}}-T_{\mathrm{prev}}$
  \State $d_k \leftarrow e^{-\eta k}\,\mathrm{tr}(T_{\mathrm{next}})$
  \State $T_{\mathrm{prev}}\!\leftarrow T_{\mathrm{curr}}$, $T_{\mathrm{curr}}\!\leftarrow T_{\mathrm{next}}$
\EndFor
\State \Return $\phi_K(A)\leftarrow [d_0,\dots,d_{K-1}]/\|[d_0,\dots,d_{K-1}]\|_2$
\end{algorithmic}
\end{algorithm}

\subsection{Adaptive Spectral Fingerprints (ASF)}
ASF produces $\phi_{K_\star}(A)$ by a \emph{dual stopping rule}:
(i) an \textbf{energy-tail} criterion; and (ii) a \textbf{Hankel low-rank} detector.

\paragraph{Energy-tail rule.}
Let $E_k=\sum_{j=0}^{k}d_j^2$ and $\rho_k=d_k^2/(E_k+\varepsilon)$ with a tiny stabilizer $\varepsilon>0$. If $\rho_k<\tau$ for $w$ consecutive $k\ge K_{\min}$, stop at $K_\star=k+1$.

\paragraph{Hankel low-rank rule (with stabilization).}
Form a Hankel $H$ from $\{d_0,\dots,d_k\}$ with $H_{ij}=d_{i+j}$, $n_H=\lfloor (k+1)/2\rfloor$. For numerical stability we (i) L2-normalize $(d_j)$ before forming $H$, (ii) optionally use a \emph{scaled Hankel} (row/column normalization), and (iii) estimate the ratio $r_H=\sigma_{\min}(H_\lambda)/\sigma_{\max}(H_\lambda)$ via Tikhonov-regularized $H_\lambda=\begin{bmatrix}H\\ \sqrt{\lambda}I\end{bmatrix}$ (e.g., $\lambda\approx 10^{-10}$) or incremental/randomized SVD. If $r_H<\varepsilon_{\mathrm{H}}$ for $w$ steps and $k\ge K_{\min}$, stop.

\paragraph{Decision logic.}
Take the earliest $k\!\ge\!K_{\min}$ that satisfies either rule for $w$ consecutive steps.

\begin{algorithm}[H]
\caption{ASF (Adaptive, exact traces)}
\label{alg:asf}
\begin{algorithmic}[1]
\Require $A$, $K_{\min}$, $K_{\max}$, $\eta$, thresholds $(\tau,\varepsilon_{\mathrm{H}})$, window $w$, $w_0$ (default $=1$)
\State $\widetilde A \leftarrow \textsc{ScaleToUnitInterval}(A)$
\Comment{for symmetric $A$; non-Hermitian uses unit-disk scaling}
\State $T_{\mathrm{prev}}\!\leftarrow I$, $T_{\mathrm{curr}}\!\leftarrow \widetilde A$
\State $d_0\!\leftarrow w_0$; $E_0\!\leftarrow d_0^2$; $h\!\leftarrow 0$
\If{$K_{\min}>1$} \State $d_1\!\leftarrow e^{-\eta}\,\mathrm{tr}(\widetilde A)$; $E_1\!\leftarrow E_0+d_1^2$ \EndIf
\For{$k=2$ to $K_{\max}-1$}
  \State $T_{\mathrm{next}}\!\leftarrow 2\widetilde A\,T_{\mathrm{curr}}-T_{\mathrm{prev}}$
  \State $d_k \leftarrow e^{-\eta k}\,\mathrm{tr}(T_{\mathrm{next}})$; $E_k\!\leftarrow E_{k-1}+d_k^2$
  \State $\rho_k \leftarrow d_k^2/(E_k+\varepsilon)$
  \State Build stabilized $H$ from $\{d_0,\dots,d_k\}$; $r_H\!\leftarrow \sigma_{\min}(H_\lambda)/\sigma_{\max}(H_\lambda)$
  \State $\textsc{hit} \leftarrow \bigl(\rho_k<\tau \,\vee\, r_H<\varepsilon_{\mathrm{H}}\bigr)\wedge (k+1\ge K_{\min})$
  \State $h \leftarrow \textsc{hit} ? (h+1) : 0$
  \If{$h \ge w$} \State $K_\star\!\leftarrow k+1$; \textbf{break} \EndIf
  \State $T_{\mathrm{prev}}\!\leftarrow T_{\mathrm{curr}}$, $T_{\mathrm{curr}}\!\leftarrow T_{\mathrm{next}}$
\EndFor
\State \Return $\phi_{K_\star}(A)\leftarrow [d_0,\dots,d_{K_\star-1}] / \|[d_0,\dots,d_{K_\star-1}]\|_2$
\end{algorithmic}
\end{algorithm}

\subsection{Hutchinson–Chebyshev for Large Sparse Matrices}
\label{sec:hutch}
For large $A$, traces are estimated with Hutchinson probes:
draw $z^{(i)}\in\{\pm 1\}^n$ (or $z^{(i)}\sim\mathcal{N}(0,I)$) and set
\[
\widehat{\mathrm{tr}}\!\big(T_k(\widetilde A)\big)=\frac{1}{p}\sum_{i=1}^p \bigl(z^{(i)}\bigr)^\top u_k^{(i)},\quad
u_0^{(i)}=z^{(i)},\; 
\]
\[
u_1^{(i)}=\widetilde A z^{(i)},\; u_{k+1}^{(i)}=2\widetilde A u_k^{(i)}-u_{k-1}^{(i)}.
\]
Define $\widehat d_k=e^{-\eta k}\widehat{\mathrm{tr}}\,T_k(\cdot)$ and the SE-guarded energy ratio
\[
\widehat{\rho}_k=\frac{\widehat d_k^2}{\sum_{j\le k}\widehat d_j^2+\varepsilon}
\quad\text{and stop if}\quad
\widehat{\rho}_k < \tau\bigl(1+\gamma \cdot \mathrm{relSE}_k \bigr),
\]
where $\mathrm{relSE}_k=\mathrm{se}_k/(|\widehat d_k|+\varepsilon)$, $\gamma\in[2,3]$.
\emph{We use a small window $w$ and a fixed $\gamma$ as a simple sequential guard; more formal error-rate control (e.g., $\alpha$-spending) is left to future work.}

\begin{algorithm}[H]
\caption{ASF-H (Adaptive, Hutchinson)}
\label{alg:asf-h}
\begin{algorithmic}[1]
\Require $A$, $p$ probes, $K_{\min}$, $K_{\max}$, $\eta$, $(\tau,\varepsilon_{\mathrm{H}})$, $w$, SE-guard $\gamma$, $w_0$ (default $=1$)
\State $\widetilde A \leftarrow \textsc{ScaleToUnitInterval}(A)$
\State For each probe $i$: $u^{(i)}_0\!\leftarrow z^{(i)}$, $u^{(i)}_1\!\leftarrow \widetilde A z^{(i)}$
\State $\widehat d_0\!\leftarrow w_0$, $\widehat E_0\!\leftarrow \widehat d_0^2$, $h\!\leftarrow 0$
\If{$K_{\min}>1$}
  \State $\widehat t_1\leftarrow \frac{1}{p}\sum_i \langle z^{(i)},u^{(i)}_1\rangle$, \quad
        $\widehat d_1\leftarrow e^{-\eta}\widehat t_1$, \quad $\widehat E_1\leftarrow \widehat E_0+\widehat d_1^2$
\EndIf
\For{$k=2$ to $K_{\max}-1$}
  \For{each $i$} \State $u^{(i)}_{k+1}\!\leftarrow 2\widetilde A u^{(i)}_{k}-u^{(i)}_{k-1}$ \EndFor
  \State $\widehat t_k \leftarrow \frac{1}{p}\sum_i \langle z^{(i)},u^{(i)}_{k}\rangle$, \quad
        $\widehat d_k \leftarrow e^{-\eta k}\widehat t_k$, \quad
        $\widehat E_k\!\leftarrow \widehat E_{k-1}+\widehat d_k^2$
  \State Estimate $\mathrm{se}_k$; set $\mathrm{relSE}_k\!\leftarrow \mathrm{se}_k/(|\widehat d_k|+\varepsilon)$
  \State Build stabilized $H$ from $\{\widehat d_0,\dots,\widehat d_k\}$; $r_H\!\leftarrow \sigma_{\min}(H_\lambda)/\sigma_{\max}(H_\lambda)$
  \State $\textsc{hit} \leftarrow \Bigl(\widehat d_k^2/\widehat E_k < \tau(1+\gamma \cdot \mathrm{relSE}_k)\;\vee\; r_H<\varepsilon_{\mathrm{H}}\Bigr)\wedge (k+1\ge K_{\min})$
  \State $h \leftarrow \textsc{hit} ? (h+1) : 0$
  \If{$h \ge w$} \State $K_\star\!\leftarrow k+1$; \textbf{break} \EndIf
\EndFor
\State \Return $\phi_{K_\star}(A)\leftarrow [\widehat d_0,\dots,\widehat d_{K_\star-1}] / \|[\widehat d_0,\dots,\widehat d_{K_\star-1}]\|_2$
\end{algorithmic}
\end{algorithm}

\begin{remark}[Scope of the claim]
Our claims concern avoiding \emph{eigenvectors or full eigen-spectra}. 
Scalar endpoint/radius estimates (via power/Lanczos, Gershgorin, or small-matrix eigenvalue calls) 
do not contradict this scope and are not required by the method.
\end{remark}

\subsection{Similarity Metrics and Multi-view Postprocessing}
By default we use Euclidean distance between $\ell_2$-normalized fingerprints; cosine is reported as an ablation.
In the E6/E6+ recommender experiments we used \emph{cosine} distance (Euclidean is included as an ablation).
When class covariances are heterogeneous, a whitened (Mahalanobis) metric can improve separation.
For nonsymmetric inputs, use $H=\tfrac12(A+A^\top)$ or $A^\top A$; multiple views (e.g., adjacency and normalized Laplacian) can be concatenated or fused by late averaging (see Appendix for details).

\subsection{Complexity and Memory}
Let $c_\mathrm{mv}$ be the cost of one sparse matvec with $A$.
CSF (exact trace): $\mathcal{O}(K\,c_\mathrm{mv})$ if traces are available; at scale we use Hutchinson.
ASF-H: $\mathcal{O}(pK\,c_\mathrm{mv})$ time and, with streaming/batching, $\mathcal{O}(n)$ memory (otherwise $\mathcal{O}(pn)$ if all probe states are kept simultaneously).
Small Hankel updates are $\mathcal{O}(K^2)$ and negligible for $K\!\le\!64$.

\section{Theory: Properties and Guarantees}
\label{sec:theory}

\paragraph{Setup.}
Let $A$ be real symmetric (or symmetrized as in Sec.~\ref{sec:prelim}) and define $\widetilde A$ by Eq.~\eqref{eq:affine} so that $\mathrm{spec}(\widetilde A)\subset[-1,1]$.
Let $T_k$ be Chebyshev polynomials, $s_k=e^{-\eta k}\,\mathrm{tr}\!\big(T_k(\widetilde A)\big)$, $E_k=\sum_{j=0}^k s_j^2$, $e_k=s_k^2/(E_k+\varepsilon)$, and let $H_k$ be the Hankel from $(s_0,\dots,s_k)$ with ratio $r_k=\sigma_{\min}(H_k)/\sigma_{\max}(H_k)$.
The fingerprint is $\phi_{K}(A)=[s_0,\dots,s_{K-1}]/\|[s_0,\dots,s_{K-1}]\|_2$.

\begin{proposition}[Invariance]
\label{prop:invariance-main}
For any invertible $S$, permutation $P$, positive diagonal $D$, and $\alpha>0$, moments (hence $\phi_K$) are invariant:
\[
s_k(S^{-1}AS)=s_k(A),\quad s_k(P^\top AP)=s_k(A),
\]
\[
s_k(D^{-1}AD)=s_k(A),\quad s_k(\alpha A)=s_k(A).
\]
\emph{Proof:} App.~\ref{app:proofs}.
\end{proposition}

\begin{theorem}[Adaptive compactness: stopping \& tail control]
\label{thm:stopping}
Fix $\tau_{\mathrm e}\!\in(0,1)$, $\tau_{\mathrm h}\!>0$, window $w$, $K_{\min}\!\ge1$.
If either (i) $e_k<\tau_{\mathrm e}$ or (ii) $r_k<\tau_{\mathrm h}$ holds for $w$ consecutive $k\!\ge\!K_{\min}$, Algorithm~\ref{alg:asf} returns $K^\star\le k^\star{+}w$ for some such $k^\star$; moreover, under (i) the normalized tail energy obeys
\[
\frac{\sum_{j>t} s_j^2}{E_t+\sum_{j>t} s_j^2}\;\le\;\frac{\tau_{\mathrm e}}{1-\tau_{\mathrm e}}.
\]
\emph{Proof:} App.~\ref{app:proofs} (uses Prop.~\ref{prop:energy}).
\end{theorem}

\begin{proposition}[Hankel low rank: mixture-of-modes]
\label{prop:hankel-lowrank}
Suppose $(s_k)$ arises from a mixture of $r$ damped exponential/cosine modes with noise variance $\sigma^2$:
$s_k=\sum_{i=1}^r a_i \rho_i^k \cos(k\theta_i+\varphi_i)+\xi_k$, with $|\rho_i|\le 1$ and $\mathbb{E}\,\xi_k=0$, $\mathrm{Var}(\xi_k)=\sigma^2$.
Then in the noiseless case $(\sigma=0)$ one has $\mathrm{rank}(H_L)\le r$ for all $L$ (under the convention that one cosine mode counts a complex-conjugate pair as a single mode). With noise, for any fixed $L$ we have $r_L=\sigma_{\min}(H_L)/\sigma_{\max}(H_L)\to 0$ as $\sigma\to 0$.
Consequently, for any $\tau_{\mathrm h}>0$ there exists $\sigma_0=\sigma_0(\tau_{\mathrm h},L,r,\{a_i,\rho_i,\theta_i,\varphi_i\})>0$ such that if $\sigma<\sigma_0$ then $r_L<\tau_{\mathrm h}$.
\emph{Proof:} App.~\ref{app:proofs}.
\end{proposition}

\begin{remark}[Counting convention]
A single real cosine mode $a\,\rho^k\cos(k\theta+\varphi)$ equals the real part of two complex exponentials with conjugate nodes.
If one counts complex-conjugate pairs separately, the Hankel rank bound becomes $\mathrm{rank}(H_L)\le 2r$; our statement uses the common convention that a conjugate pair is counted as one ``mode'', hence the bound $\le r$.
\end{remark}

\begin{theorem}[Hutchinson concentration]
\label{thm:hutch-conc}
Let $B_k=T_k(\widetilde A)$ and $\widehat{t}_k=\tfrac{1}{p}\sum_{i=1}^p (z^{(i)})^\top B_k z^{(i)}$ with i.i.d.\ Rademacher probes. For any $\delta\in(0,1)$,
\[
\bigl|\widehat{t}_k-\mathrm{tr}(B_k)\bigr|\le \|B_k\|_F\sqrt{\tfrac{2\log(2/\delta)}{p}},
\]
\[
\bigl|\widehat{s}_k-s_k\bigr|\le e^{-\eta k}\|B_k\|_F\sqrt{\tfrac{2\log(2/\delta)}{p}}.
\]
(Uniform-in-$k$ version up to $K$ in App., Cor.~\ref{cor:uniform-k}.) \emph{Proof:} App.~\ref{app:proofs}.
\end{theorem}

\begin{theorem}[Lipschitz stability]
\label{thm:lipschitz}
If $\mathrm{spec}(\widetilde A),\mathrm{spec}(\widetilde B)\subset[-1,1]$, then for all $k\ge0$,
\[
\bigl|\mathrm{tr}\,T_k(\widetilde A)-\mathrm{tr}\,T_k(\widetilde B)\bigr|
\le n\,C_T\,k^2\,\|\widetilde A-\widetilde B\|_2\quad(C_T\le2),
\]
and for any $K\ge1$,
\[
\bigl\|\phi_{K}(A)-\phi_{K}(B)\bigr\|_2
\le 
\]
\[
\frac{n}{\|[s_0,\dots,s_{K-1}]\|_2}\Bigl(\sum_{k=0}^{K-1} e^{-2\eta k} C_T^2 k^4 \Bigr)^{\!1/2}\!\|\widetilde A-\widetilde B\|_2
\]
\[
= \mathcal{O}\!\big(n\,\eta^{-5/2}\,\|\widetilde A-\widetilde B\|_2\big).
\]
\emph{Proof:} App.~\ref{app:proofs}.
\end{theorem}

\section{Experiments}
\label{sec:experiments}

We evaluate on synthetic suites (E0--E2), strong baselines (E2b), a SuiteSparse mini-benchmark (E3), Hutchinson ablations (E4), noise stability (E5), and an adversarial preconditioner-selection stress test (E6/E6+). Unless noted, we use Euclidean distance on $\ell_2$-normalized fingerprints (cosine is an ablation), CSF with $w_0{=}1$ and small damping $\eta>0$.

\subsection{E0: Invariance \& Scaling (tables only)}
We test permutation, (diagonal/general) similarity, and positive scaling with/without our spectral normalization. With normalization (\emph{scaled}), all transforms collapse to $\approx 0$ distance (below the numerical floor); without normalization, only positive scaling (\texttt{alpha}) breaks invariance.

\begin{table}[H]
\centering\small
\caption{E0 (scaled): distance to original after transform (mean/median/IQR). All are $\approx 0$ (below numerical floor, $<10^{-15}$).}
\begin{tabular}{lccc}\toprule
transform & mean & median & IQR\\\midrule
alpha & 0.0 & 0.0 & 0.0\\
diag\_sim & 0.0 & 0.0 & 0.0\\
gen\_sim & 0.0 & 0.0 & 0.0\\
perm & 0.0 & 0.0 & 0.0\\\bottomrule
\end{tabular}
\end{table}

\begin{table}[H]
\centering\small
\caption{E0 (no-scale): positive scaling breaks invariance; others remain zero.}
\begin{tabular}{lccc}\toprule
transform & mean & median & IQR\\\midrule
alpha & 0.9782 & 1.3316 & 1.3015\\
diag\_sim & 0.0000 & 0.0000 & 0.0000\\
gen\_sim & 0.0000 & 0.0000 & 0.0000\\
perm & 0.0000 & 0.0000 & 0.0000\\\bottomrule
\end{tabular}
\end{table}

\subsection{E1: 4-Family Compactness (tables only)}
We sweep $K$ and compare CSF/ASF to spectral/entrywise baselines. CSF-$K{=}3,5$ achieve \textbf{ARI $=1.0$} with top silhouettes and lowest runtimes; ASF (exact) also reaches ARI $=1.0$; Hutchinson (ASF-H) preserves accuracy at higher cost.

\begin{table}[H]
\centering\small
\caption{E1 K-sweep: ARI / silhouette / runtime (s).}
\resizebox{\linewidth}{!}{%
\begin{tabular}{lcccc}\toprule
method & $K$ or $m$ & ARI & Silhouette & runtime\\\midrule
ASF (Adaptive, exact trace) & -- & 1.0000 & 0.8510 & 0.0692\\
ASF-H (Adaptive, Hutchinson) & -- & 1.0000 & 0.8555 & 4.2217\\
CSF-K (Fixed-K Chebyshev) & 3 & 1.0000 & 0.8928 & 0.0431\\
CSF-K (Fixed-K Chebyshev) & 5 & 1.0000 & 0.8942 & 0.0457\\
CSF-K (Fixed-K Chebyshev) & 10 & 1.0000 & 0.8550 & 0.0632\\
CSF-K (Fixed-K Chebyshev) & 50 & 1.0000 & 0.8100 & 0.2068\\
Baseline: Raw Power Moments & 10 & 1.0000 & 0.8755 & 0.0562\\
Baseline: Spectral Norm $\|\cdot\|_2$ & -- & 0.6508 & 0.7505 & --\\
Baseline: Top-$m$ Eigenvalues & 16 & 0.6220 & 0.6814 & 0.0635\\
Baseline: Heat-Trace @ $T$ & -- & 0.5041 & 0.5492 & 0.0584\\
Baseline: Frobenius & -- & 0.3260 & 0.3028 & --\\
CSF-K (Fixed-K Chebyshev) & 1 & 0.0020 & 0.0000 & 0.0598\\\bottomrule
\end{tabular}}
\end{table}

\begin{table}[H]
\centering\small
\caption{E1: Adaptive $K^\star$ summary (ASF).}
\begin{tabular}{lccc}\toprule
 & mean & median & IQR\\\midrule
$K^\star$ & 8.39 & 9.0 & 3.0\\\bottomrule
\end{tabular}
\end{table}

\subsection{E2: 5-Family (BA vs.\ ER included)}
On a harder 5-family suite, \textbf{CSF-$K{=}5$} and \textbf{ASF} both reach \textbf{ARI $=1.0$}. We report overall and per-view silhouettes.

\begin{table}[H]
\centering
\caption{E2: Per-method scores.}
\resizebox{\linewidth}{!}{%
\begin{tabular}{lccccccc}\toprule
method & ARI & Silh. (overall) & Silh.[Cov] & Silh.[Kernel] & Silh.[GOE] & Silh.[Adj-BA] & Silh.[Adj-ER]\\\midrule
ASF (Adaptive, exact) & 1.0000 & 0.7989 & 0.8160 & 0.8207 & 0.8591 & 0.7248 & 0.7737\\
CSF-K (Fixed-K) & 1.0000 & 0.8209 & 0.9054 & 0.8642 & 0.8466 & 0.7699 & 0.7184\\
Baseline: Top-$m$ Eigen & 0.6933 & 0.5855 & 0.8850 & 0.8126 & 0.6057 & 0.2727 & 0.3513\\
Baseline: Heat-Trace & 0.4583 & 0.6048 & 0.8502 & 0.7710 & 0.4504 & 0.4569 & 0.4957\\\bottomrule
\end{tabular}}
\end{table}

\begin{figure}[H]
\centering
\includegraphics[width=\linewidth]{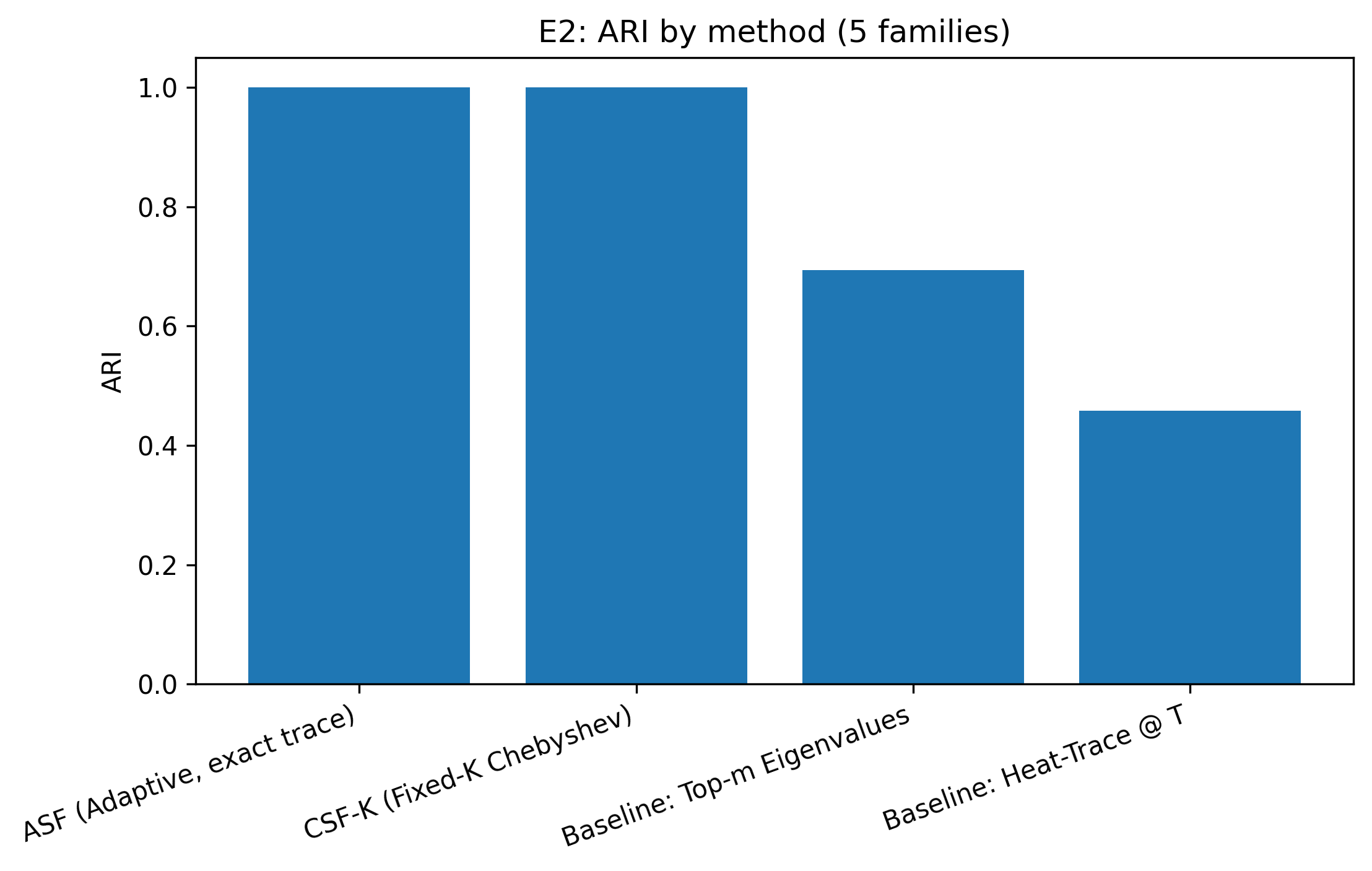}
\caption{E2: ARI/silhouette overview across methods and views.}
\end{figure}

\begin{figure}[H]
\centering
\includegraphics[width=\linewidth]{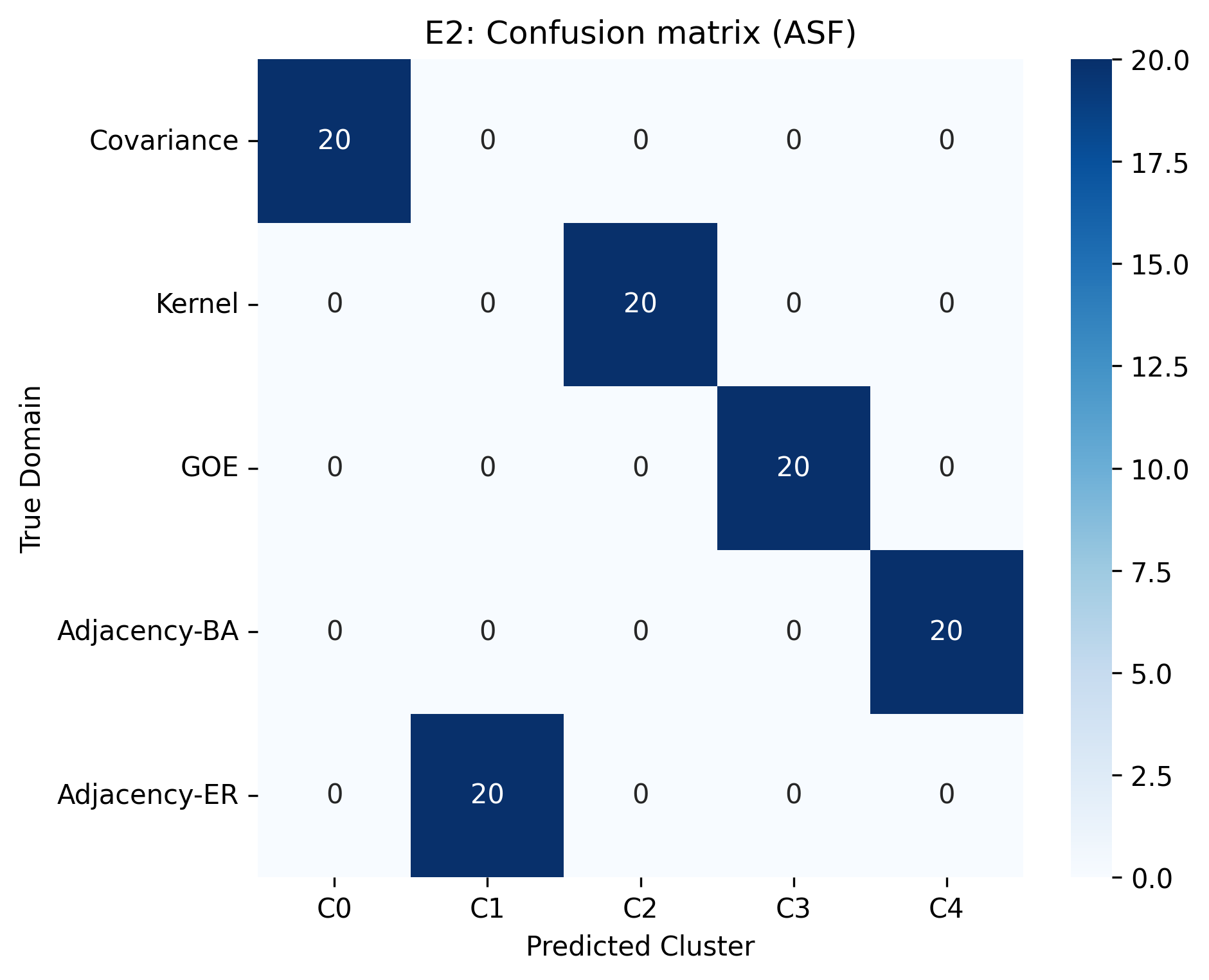}
\caption{E2: Confusion matrix (one run; label permutations allowed). ARI $=1.0$.}
\end{figure}

\subsection{E2b: Strong Baselines}
CSF-$K{=}5$ matches the best accuracy (ARI $=1.0$) while using only 5 dimensions and avoiding eigendecomposition.

\begin{table}[H]
\centering\small
\caption{E2b: Accuracy vs.\ feature dimension vs.\ time; eigendecomposition required?}
\resizebox{\linewidth}{!}{%
\begin{tabular}{lccccc}\toprule
method & ARI & Silh. & Dim & runtime (s) & Eigendecomp\\\midrule
CSF-K=5 (ours) & 1.000 & 0.8209 & 5 & 0.0604 & No\\
ASF (ours) & 1.000 & 0.7989 & 8 & 0.0875 & No\\
EigenHist+Wasserstein & 1.000 & 0.7738 & 64 & 0.3573 & Yes\\
WL-Subtree (graphs) & 1.000 & 0.3067 & 9153 & 0.0582 & No\\
HeatKernel (L2 on traces) & 0.773 & 0.7005 & 5 & 0.0758 & Yes\\\bottomrule
\end{tabular}}
\end{table}

\subsection{E3: SuiteSparse Mini-Benchmark (Hutchinson)}
On four real matrices, both CSF-H and ASF-H achieve \textbf{ARI $=1.0$}, validating sketch-based fingerprints. Specifically, we use \texttt{HB/bcsstk01}, \texttt{HB/bcsstk06}, \texttt{HB/gr\_30\_30}, and \texttt{AG-Monien/netz4504} from the SuiteSparse Matrix Collection.

\begin{table}[H]
\centering\small
\caption{E3: SuiteSparse mini-benchmark.}
\begin{tabular}{lcc}\toprule
method & ARI & Silhouette\\\midrule
CSF-K (Hutchinson) & 1.0 & 0.6220\\
ASF-H (Adaptive, Hutchinson) & 1.0 & 0.5386\\\bottomrule
\end{tabular}
\end{table}

\begin{figure}[H]
\centering
\includegraphics[width=\linewidth]{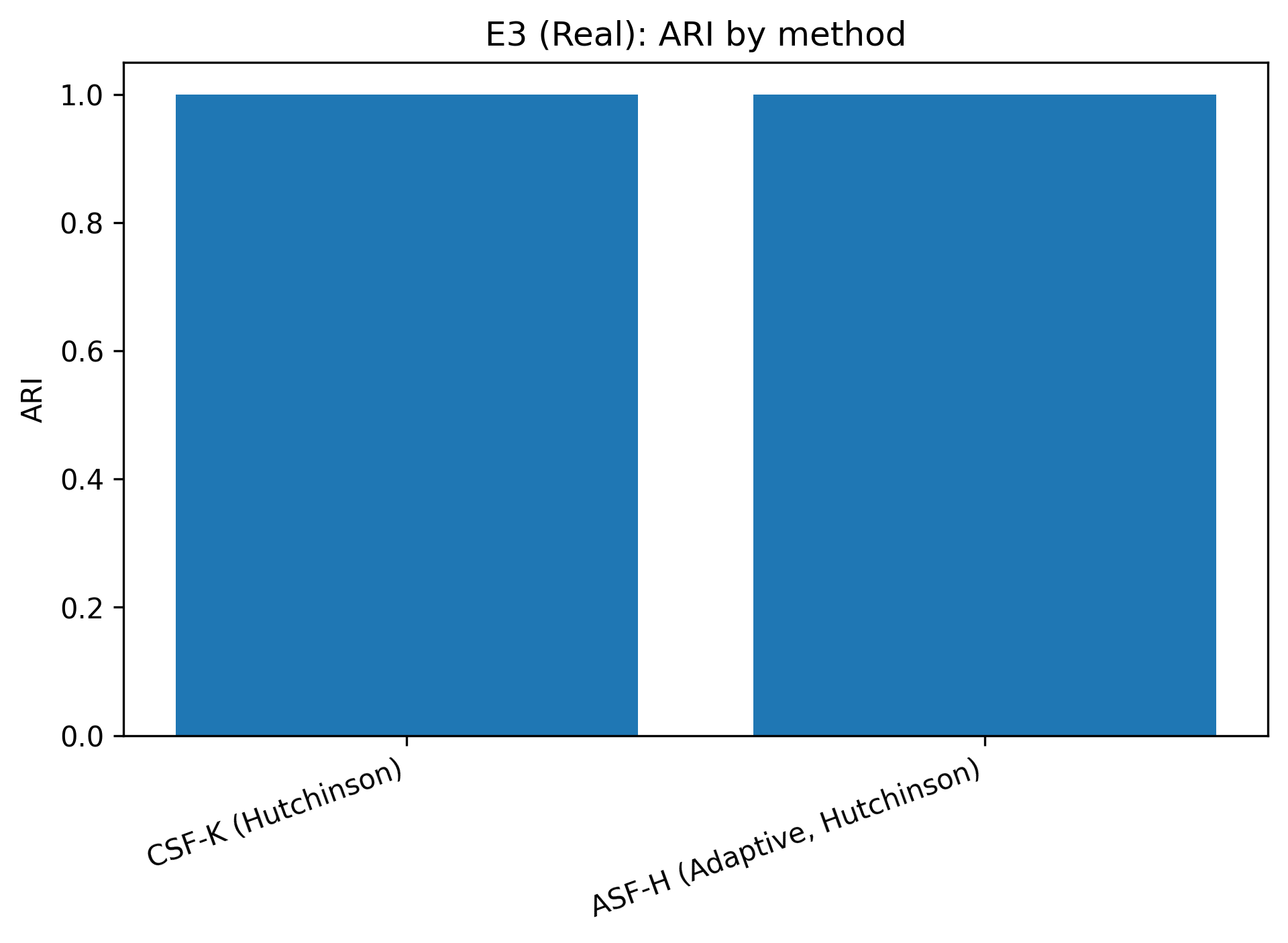}
\caption{E3: ARI on SuiteSparse (Hutchinson).}
\end{figure}

\subsection{E4: Hutchinson Ablation (probe count $p$)}
ARI remains at \textbf{1.0} across all $p$ and $K\!\in\!\{3,5\}$; silhouette increases smoothly; runtime scales near-linearly.

\begin{table}[H]
\centering\small
\caption{E4: ARI vs.\ $p$ (columns are $K$).}
\begin{tabular}{lcc}\toprule
$p$ & $K{=}3$ & $K{=}5$\\\midrule
10 & 1.0 & 1.0\\
50 & 1.0 & 1.0\\
100 & 1.0 & 1.0\\
200 & 1.0 & 1.0\\
500 & 1.0 & 1.0\\
1000 & 1.0 & 1.0\\\bottomrule
\end{tabular}
\end{table}

\begin{table}[H]
\centering\small
\caption{E4: Silhouette vs.\ $p$ (columns are $K$).}
\begin{tabular}{lcc}\toprule
$p$ & $K{=}3$ & $K{=}5$\\\midrule
10 & 0.8740 & 0.8587\\
50 & 0.8894 & 0.8801\\
100 & 0.8899 & 0.8899\\
200 & 0.8910 & 0.8904\\
500 & 0.8914 & 0.8946\\
1000 & 0.8935 & 0.8941\\\bottomrule
\end{tabular}
\end{table}

\begin{table}[H]
\centering\small
\caption{E4: Runtime (s) vs.\ $p$ (columns are $K$).}
\begin{tabular}{lcc}\toprule
$p$ & $K{=}3$ & $K{=}5$\\\midrule
10 & 0.0787 & 0.0893\\
50 & 0.1896 & 0.2410\\
100 & 0.3223 & 0.4385\\
200 & 0.5984 & 0.8278\\
500 & 1.3878 & 1.9215\\
1000 & 2.7560 & 2.7219\\\bottomrule
\end{tabular}
\end{table}

\begin{figure}[H]
\centering
\includegraphics[width=\linewidth]{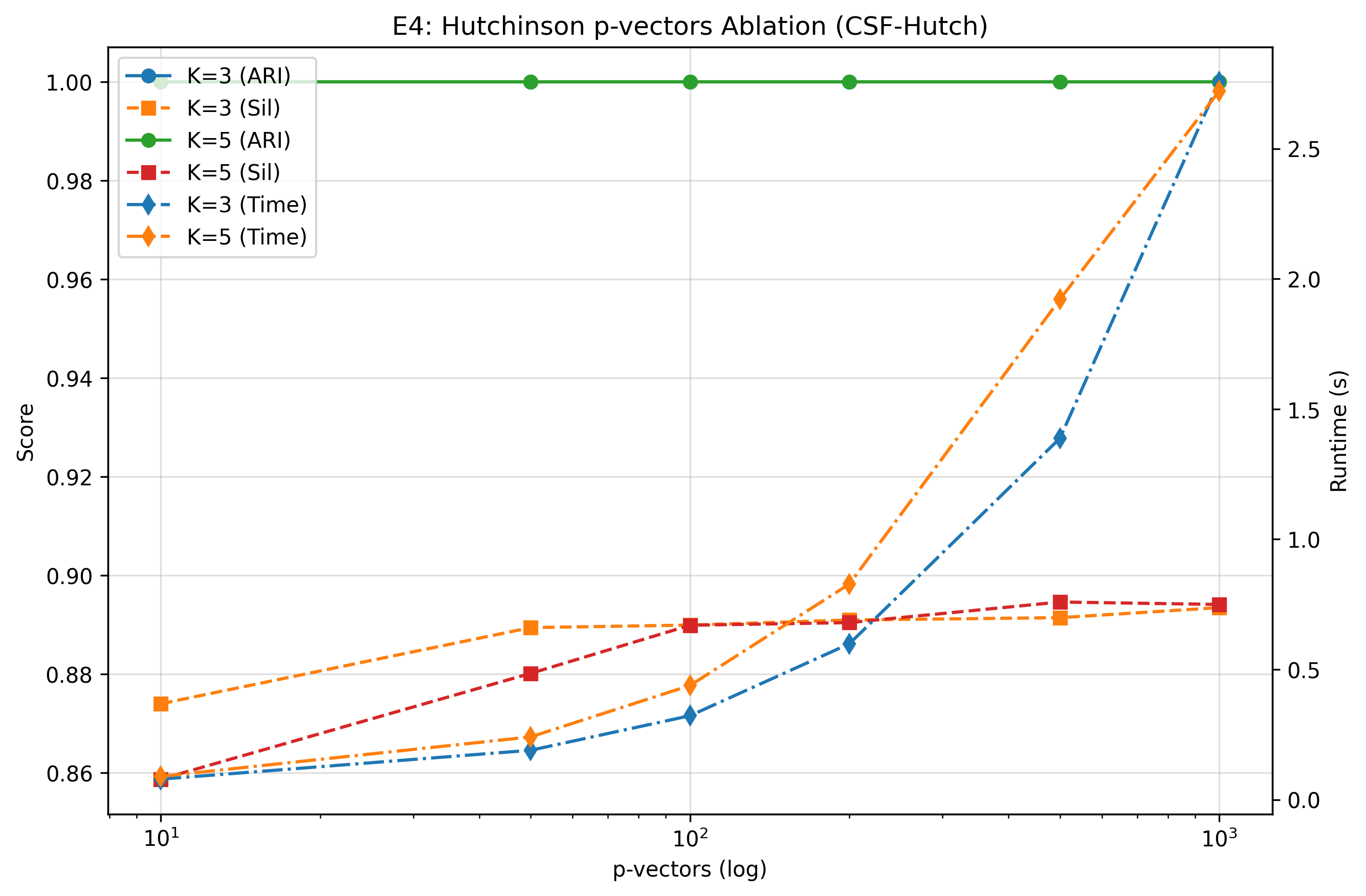}
\caption{E4: Probe--quality--time trade-off for CSF-H/ASF-H.}
\end{figure}

\subsection{E5: Noise Stability}
Fingerprint distance scales $\propto \epsilon^{1.03}$ with additive noise; fit $R^2 \approx 0.993$.

\begin{table}[H]
\centering\small
\caption{E5: Noise--distance regression (log--log).}
\begin{tabular}{lc}\toprule
Metric & Value\\\midrule
Slope & 1.0289\\
$R^2$ & 0.9931\\\bottomrule
\end{tabular}
\end{table}

\begin{figure}[H]
\centering
\includegraphics[width=\linewidth]{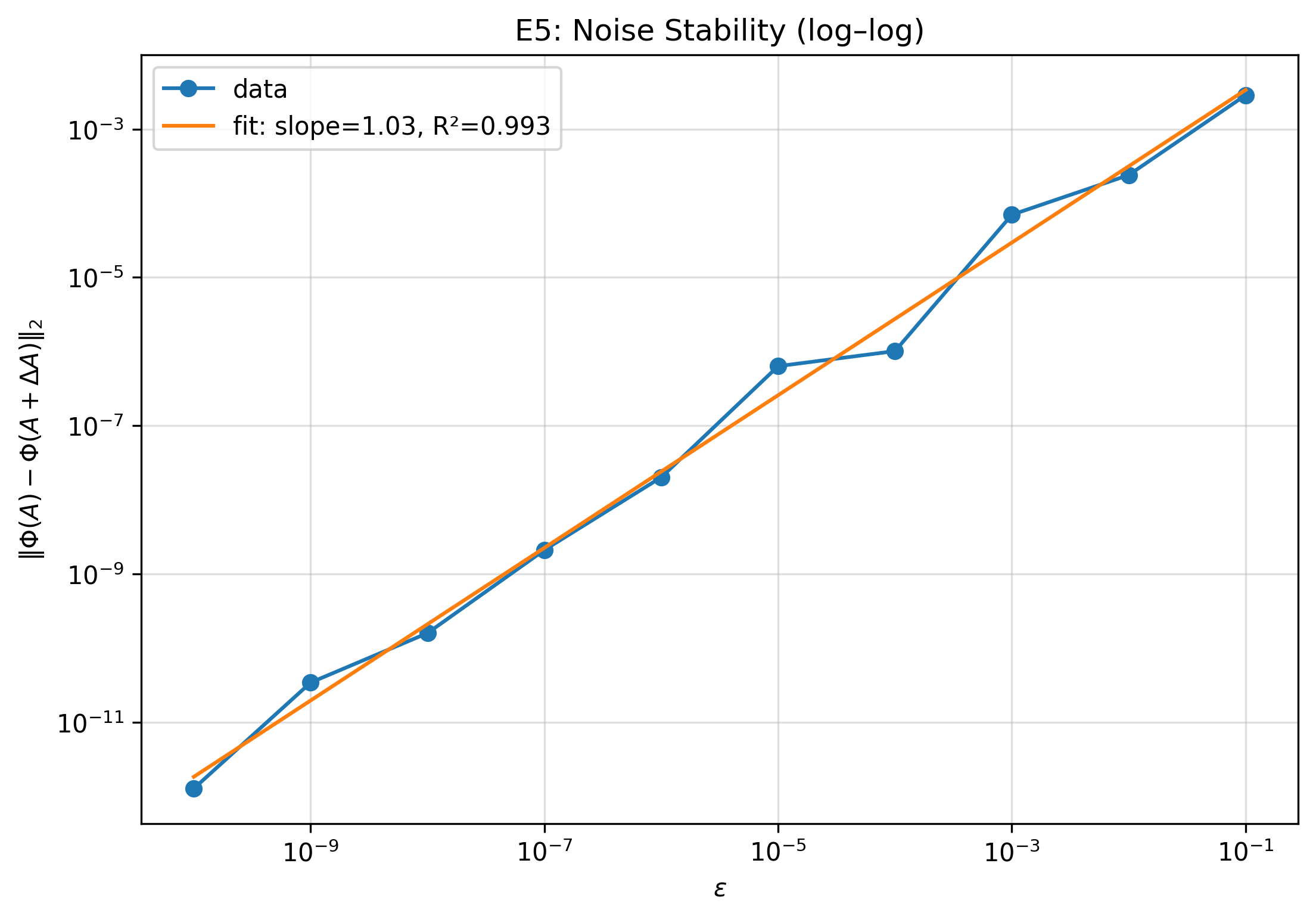}
\caption{E5: Log--log fit of fingerprint distance vs.\ noise level $\epsilon$.}
\end{figure}

\subsection{E6/E6+: Trap-Robust Preconditioner Selection}
Across adversarial seeds, the phylogeny-guided policy attains \textbf{success rate $1.00$} with \textbf{p90 regret $=0$}; Frobenius 1-NN exhibits \textbf{p90 $=34$} extra iterations.

\begin{table}[H]
\centering\small
\caption{E6+: Aggregated adversarial results (success and regret).}
\begin{tabular}{lc}\toprule
Metric & Value\\\midrule
Phy Success Rate & 1.00\\
Fro-1NN Success Rate & 0.75\\
Fro-$k$NN Success Rate & 0.95\\
Phy Extra Iters (Median) & 0.00\\
Fro-1NN Extra Iters (Median) & 0.00\\
Fro-$k$NN Extra Iters (Median) & 0.00\\
Phy Extra Iters (Mean) & 0.00\\
Fro-1NN Extra Iters (Mean) & 5.70\\
Fro-$k$NN Extra Iters (Mean) & 4.55\\
Phy Extra Iters (p90) & 0.00\\
Fro-1NN Extra Iters (p90) & 34.00\\
Fro-$k$NN Extra Iters (p90) & 0.00\\\bottomrule
\end{tabular}
\end{table}

\begin{figure}[H]
\centering
\includegraphics[width=\linewidth]{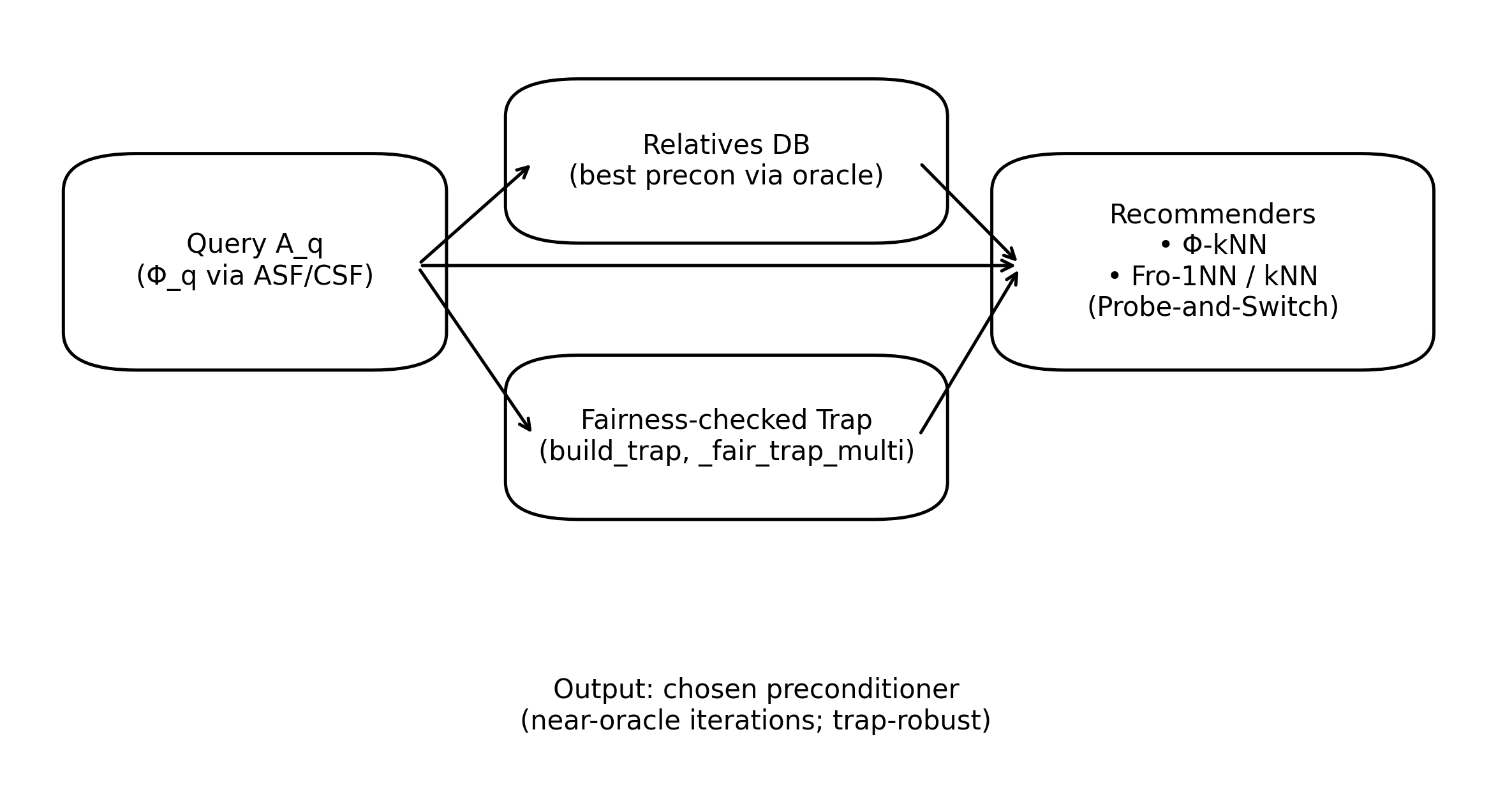}
\caption{E6: \emph{Probe-and-Switch} flow (trap detection and switching).}
\end{figure}


\section{Related Work}

\paragraph{Spectral fingerprints and graph descriptors.}
Permutation-invariant spectral summaries are longstanding. Heat-trace descriptors capture multiscale diffusion (e.g., NetLSD) and compare graphs without node alignment \citep{Tsitsulin2018}. Eigenvalue histograms compared by optimal transport (e.g., Wasserstein) are effective but typically require (partial) eigendecompositions and careful binning \citep{Villani2008}. Weisfeiler--Lehman subtree kernels \citep{Shervashidze2011} encode rich combinatorics yet are not similarity-invariant and can be fragile outside the WL regime. Our CSF/ASF fingerprints avoid computing eigenvectors or full eigen-spectra, operate via matvecs, and are invariant to permutation, similarity, and positive scaling.

\paragraph{Chebyshev expansions and stochastic traces.}
Chebyshev/KPM expansions approximate spectral densities and traces from matrix–vector products \citep{Weisse2006}. Hutchinson’s estimator and refinements (Hutch++) enable scalable trace estimation with variance $\mathcal{O}(1/p)$ for $p$ probe vectors \citep{Hutchinson1989,Meyer2021HutchPP}. Stochastic Lanczos quadrature provides complementary density-of-states estimates \citep{Ubaru2017}. These techniques sit within randomized numerical linear algebra \citep{Halko2011}. We combine affine spectrum mapping, damped Chebyshev moments, and Hutchinson sketches into compact, similarity-invariant fingerprints.

\paragraph{Preconditioning and automated choices.}
Classical work shows the centrality of preconditioning for Krylov solvers (CG, GMRES) \citep{Hestenes1952,Meijerink1977,Trottenberg2001,Saad2003,Benzi2002}. Prior auto-tuning often relies on heavy spectral surrogates or metadata. Our \emph{matrix phylogeny} is orthogonal: a retrieval layer built from compact fingerprints that steers a probe-and-switch policy, achieving near-oracle iteration counts without eigenspectra.

\paragraph{Limits of spectral identifiability.}
Cospectral nonisomorphic graphs (and matrices) exist \citep{vanDam2003}; descriptors depending only on spectra, including ours, cannot distinguish them in principle. We mitigate empirically via multi-view fingerprints (adjacency/Laplacian/GOE/kernel) and find strong separation across families.

\section{Limitations}

\textbf{Cospectral ambiguity.}
Spectral equivalence remains a worst-case confound; multi-view concatenation helps but cannot resolve all cases.

\textbf{Endpoint/scaling estimates.}
Chebyshev stability hinges on mapping spectra to $[-1,1]$. Overly loose bounds can overdamp higher moments, sometimes requiring a slightly larger $K$ or ASF’s adaptive stop.

\textbf{Sketch variance.}
Hutchinson estimates concentrate as $1/\sqrt{p}$; high-variance instances may need more probe vectors to stabilize ASF’s stopping rule.

\textbf{Nonsymmetric inputs.}
For non-Hermitian $A$ we use radius scaling or symmetrization ($A^\top A$ or $(A{+}A^\top)/2$). This preserves invariances but may blur non-normal features when they matter.

\textbf{Actionability scope.}
Our auto-selection uses a standard preconditioner portfolio; specialized domains may benefit from broader candidates or lightweight non-spectral cues (e.g., sparsity sketches).

\section{Conclusion}

We presented \emph{Matrix Phylogeny} via \textbf{Compact/Adaptive Spectral Fingerprints} (CSF/ASF): tiny ($K{\le}10$), eigendecomposition-free, and invariant descriptors built from Chebyshev moments and Hutchinson sketches. Across E0–E6+, we observe \emph{phylogenetic compactness}: \textbf{CSF-$K{=}3$–$5$} already yields perfect clustering on synthetic and real suites while being fastest, and \textbf{ASF} provides a conservative safety net (median $K^\star{\approx}9$). The fingerprints are stable to perturbations (noise–distance slope $\approx\!1.03$, $R^2\!\approx\!0.993$) and power an actionable pipeline: a nearest-neighbor \emph{probe-and-switch} recommender attains near-oracle iterations in adversarial settings, avoiding the trap behavior seen with Frobenius baselines.

The recipe—affine scaling to $[-1,1]$, Chebyshev recurrences, Hutchinson sketches, cosine distances—unlocks scalability, invariance, and accuracy with feature dimensions far smaller than spectral baselines. We recommend \textbf{CSF-$K{=}5$} as the default and \textbf{ASF} when domain adaptivity is required. Future work: principled multi-view fusion, solver-in-the-loop tuning of damping/stopping/distances for regret, and identifiability theory for when compact moments uniquely determine matrix families.


\newpage

\bibliographystyle{apalike}
\bibliography{references}

\section*{Reproducibility Checklist}

\begin{enumerate}

  \item For all models and algorithms presented, check if you include:
  \begin{enumerate}
    \item A clear description of the mathematical setting, assumptions, algorithm, and/or model. [Yes]
    \item An analysis of the properties and complexity (time, space, sample size) of any algorithm. [Yes]
    \item (Optional) Anonymized source code, with specification of all dependencies, including external libraries. [Yes]
  \end{enumerate}

  \item For any theoretical claim, check if you include:
  \begin{enumerate}
    \item Statements of the full set of assumptions of all theoretical results. [Yes]
    \item Complete proofs of all theoretical results. [Yes]
    \item Clear explanations of any assumptions. [Yes]     
  \end{enumerate}

  \item For all figures and tables that present empirical results, check if you include:
  \begin{enumerate}
    \item The code, data, and instructions needed to reproduce the main experimental results (either in the supplemental material or as a URL). [Yes]
    \item All the training details (e.g., data splits, hyperparameters, how they were chosen). [Yes]
    \item A clear definition of the specific measure or statistics and error bars (e.g., with respect to the random seed after running experiments multiple times). [Yes]
    \item A description of the computing infrastructure used. (e.g., type of GPUs, internal cluster, or cloud provider). [Yes]
  \end{enumerate}

  \item If you are using existing assets (e.g., code, data, models) or curating/releasing new assets, check if you include:
  \begin{enumerate}
    \item Citations of the creator If your work uses existing assets. [Yes]
    \item The license information of the assets, if applicable. [Not Applicable]
    \item New assets either in the supplemental material or as a URL, if applicable. [Not Applicable]
    \item Information about consent from data providers/curators. [Not Applicable]
    \item Discussion of sensible content if applicable, e.g., personally identifiable information or offensive content. [Not Applicable]
  \end{enumerate}

  \item If you used crowdsourcing or conducted research with human subjects, check if you include:
  \begin{enumerate}
    \item The full text of instructions given to participants and screenshots. [Not Applicable]
    \item Descriptions of potential participant risks, with links to Institutional Review Board (IRB) approvals if applicable. [Not Applicable]
    \item The estimated hourly wage paid to participants and the total amount spent on participant compensation. [Not Applicable]
  \end{enumerate}

\end{enumerate}

\newpage

\clearpage
\appendix
\onecolumn

\section{Appendix: Theoretical Guarantees and Implementation Details}
\label{app:theory-impl}

\subsection{Notation and Theorem Environments}
We assume throughout that $A\in\mathbb{R}^{n\times n}$ is real symmetric (or is symmetrized as $H=\tfrac{1}{2}(A+A^\top)$ or $A^\top A$).
Let the affine spectral normalization be
\[
\widetilde A=\frac{A-mI}{r},\qquad
m=\tfrac12\bigl(\lambda_{\max}(A)+\lambda_{\min}(A)\bigr),
\]
\[
r_0=\tfrac12\bigl(\lambda_{\max}(A)-\lambda_{\min}(A)\bigr),\;\;
r=(1+\varepsilon_{\mathrm{rel}})\,r_0,
\]
with a tiny relative margin $\varepsilon_{\mathrm{rel}}>0$, so that $\mathrm{spec}(\widetilde A)\subset[-1,1]$ and scaling $A\mapsto \alpha A$ leaves $\widetilde A$ unchanged.
We write $T_k$ for the Chebyshev polynomials of the first kind, $T_0(x)=1$, $T_1(x)=x$, $T_{k+1}(x)=2xT_k(x)-T_{k-1}(x)$.
Define the (optionally damped) moments $s_k(A)=e^{-\eta k}\,\mathrm{tr}\!\big(T_k(\widetilde A)\big)$, cumulative energy $E_K(A)=\sum_{k=0}^{K}s_k(A)^2$, energy share $e_k=s_k^2/(E_k+\varepsilon)$, and Hankel ratio $r_k=\sigma_{\min}(H_k)/\sigma_{\max}(H_k)$ for $H_k=[s_{i+j}]_{i,j=0}^{\lfloor k/2\rfloor}$.
The fingerprint is $\phi_K(A)=[s_0,\dots,s_{K-1}]/\|[s_0,\dots,s_{K-1}]\|_2$, and $K^\star$ is selected by Algorithm~\ref{alg:asf}.

\makeatletter
\@ifundefined{theorem}{\newtheorem{theorem}{Theorem}}{}
\@ifundefined{proposition}{\newtheorem{proposition}{Proposition}}{}
\@ifundefined{lemma}{\newtheorem{lemma}{Lemma}}{}
\makeatother

\begin{proposition}[Invariance (appendix)]
\label{prop:invariance-app}
Let $S$ be invertible, $P$ a permutation (orthogonal), $D$ positive diagonal, and $\alpha>0$.
With the relative-margin normalization above, for all $k\ge0$,
\[
\mathrm{tr}\,T_k\!\big((S^{-1}AS)^{\;\widetilde{}}\big)
=\mathrm{tr}\,T_k(\widetilde A),
\]
\[
\mathrm{tr}\,T_k\!\big((P^\top AP)^{\;\widetilde{}}\big) =\mathrm{tr}\,T_k(\widetilde A),
\]
\[
\mathrm{tr}\,T_k\!\big((D^{-1}AD)^{\;\widetilde{}}\big)
=\mathrm{tr}\,T_k(\widetilde A),
\]
and $\mathrm{tr}\,T_k\!\big((\alpha A)^{\;\widetilde{}}\big)=\mathrm{tr}\,T_k(\widetilde A)$.
Consequently $\phi_K(\cdot)$ computed on $\widetilde A$ is invariant to similarity, permutations, diagonal similarities, and global positive scaling.
\end{proposition}

\paragraph{Proof sketch.}
Polynomials conjugate: $T_k(S^{-1}AS)=S^{-1}T_k(A)S$ and $\mathrm{tr}$ is conjugation-invariant; permutation is the orthogonal case, diagonal similarity is identical.
With $r=(1+\varepsilon_{\mathrm{rel}})r_0$ and $m,r_0$ scaling by $\alpha$, normalization is unchanged under $A\mapsto\alpha A$. \qed

\begin{lemma}[Hutchinson variance and SE-guard]
\label{lem:hutch}
Let $B_k=T_k(\widetilde A)$ and $z\in\{\pm1\}^n$ be i.i.d.\ Rademacher.
Then $\mathbb{E}[z^\top B_k z]=\mathrm{tr}(B_k)$ and $\mathrm{Var}[z^\top B_k z]\le 2\|B_k\|_F^2$.
With $p$ probes, the standard error satisfies
$\mathrm{se}_k \le \sqrt{2}\,\|B_k\|_F/\sqrt{p}$ (scaled by $e^{-\eta k}$ for $d_k$).
Consequently, the SE-guarded stopping rule
$\widehat e_k<\tau_{\mathrm{e}}\bigl(1+\gamma\,\mathrm{relSE}_k\bigr)$
reduces false stops as $p$ increases.
\end{lemma}

\begin{proposition}[Energy-tail guarantee]\label{prop:energy}
Let $d_k=e^{-\eta k}\mathrm{tr}\,T_k(\widetilde A)$ and $E_k=\sum_{j\le k} d_j^2$.
If $d_k^2/E_k<\tau$ holds for $w$ consecutive steps from some index $t$ onward
and $(d_k^2)$ is eventually nonincreasing, then
\[
\frac{\sum_{j>t} d_j^2}{E_t+\sum_{j>t} d_j^2} \;\le\; \frac{\tau}{1-\tau}.
\]
\end{proposition}

\paragraph{Rank-sensitive corollary.}
Since $\|B_k\|_F\le \sqrt{\mathrm{rank}(B_k)}\,\|B_k\|_2$ and $\|B_k\|_2\le 1$ on $[-1,1]$, one gets $p=\mathcal{O}\!\big(\mathrm{rank}(B_k)\,\varepsilon^{-2}\log(1/\delta)\big)$.

\begin{corollary}[Uniform-in-$k$ up to $K$]
\label{cor:uniform-k}
Fix $K\ge1$ and $\delta\in(0,1)$. With probability at least $1-\delta$,
for all $k=0,1,\dots,K-1$,
\[
\bigl|\widehat t_k-\mathrm{tr}(B_k)\bigr|
\;\le\;
\|B_k\|_F \sqrt{\tfrac{2\log(2K/\delta)}{p}},
\qquad
\bigl|\widehat s_k-s_k\bigr|
\;\le\;
e^{-\eta k}\|B_k\|_F \sqrt{\tfrac{2\log(2K/\delta)}{p}}.
\]
\emph{Proof.} Apply Theorem~\ref{thm:hutch-conc} with $\delta\leftarrow \delta/K$ to each $k$
and union bound over $k=0,\dots,K-1$.
\qed
\end{corollary}

\subsection{Implementation Details}
\paragraph{Spectral interval estimation.}
We estimate $(\lambda_{\min},\lambda_{\max})$ via a few (3--5) steps of power and shifted-power iterations, optionally refined by Gershgorin disks; a small margin $\delta>0$ is added.
If the gap estimate is unstable, we fall back to a conservative spectral-radius bound.

\paragraph{Endpoint fallback (relative margin).}
If $(\lambda_{\min},\lambda_{\max})$ estimation is unstable, we combine Gershgorin bounds with a few steps of (shifted) power iterations and then apply a \emph{relative} safety margin $r=(1+\varepsilon_{\mathrm{rel}})r_0$. This keeps the scaling invariance intact even under $A\mapsto \alpha A$.

\paragraph{Stable numerics.}
We add $\varepsilon\in[10^{-12},10^{-10}]$ to denominators, cap condition numbers when forming small Hankels, and use randomized SVD if dense SVD is ill-conditioned.
Chebyshev recurrences use the stable three-term update with periodic re-orthogonalization in mixed precision if needed.

\paragraph{Streaming memory.}
ASF-H stores only $(u^{(i)}_{k-1},u^{(i)}_{k})$ per probe and streams probes in micro-batches, giving memory $\mathcal{O}(pn)$ (or $\mathcal{O}(n)$ with batching).

\paragraph{Stopping defaults and knobs.}
Recommended: $\eta=0.06$, $K_{\min}=2\text{--}3$, $K_{\max}=64$, window $w=2\text{--}3$, $\tau_{\mathrm{e}}=10^{-3}$ (retrieval) or $5\!\times\!10^{-4}$ (fine discrimination), Hankel threshold $\tau_{\mathrm{h}}=10^{-3}$, Hutchinson $p=64$ with Rademacher probes, SE-guard factor $\gamma=2$.

\paragraph{Distance choices.}
Cosine on L2-normalized fingerprints is default.
For heteroskedastic classes, a class-covariance-whitened (Mahalanobis) distance can improve separation; we also report Euclidean for ablations.

\paragraph{Views of a matrix.}
When $A$ is not symmetric, use $H=\tfrac{1}{2}(A+A^\top)$ or $A^\top A$.
Multiple views (e.g., adjacency, normalized Laplacian, kernel, covariance) can be concatenated or fused by late averaging.

\subsection{Complexity Summary}
Let $\mathrm{Tmv}$ denote the cost of a matrix--vector multiply with $A$.
\begin{itemize}
\item \textbf{CSF-K (exact traces).} Cost $\mathcal{O}(K\cdot \mathrm{Tmv})$ for the Chebyshev recurrence plus one trace per order (via explicit trace only if diagonal known; otherwise use Hutchinson).
\item \textbf{CSF-K (Hutchinson).} Cost $\mathcal{O}(p\,K\cdot \mathrm{Tmv})$; trivially parallel across probes.
\item \textbf{ASF (adaptive, exact).} As above with $K$ replaced by $K^\star\le K_{\max}$ determined by the energy/Hankel stop.
\item \textbf{ASF-H (adaptive, Hutchinson).} $\mathcal{O}(p\,K^\star\cdot \mathrm{Tmv})$, plus $\mathcal{O}(K^\star)$ for Hankel/SVD on tiny $(\le \lceil K^\star/2\rceil)$ matrices.
\end{itemize}

\subsection{Limitations and Extensions}
Cospectral nonisomorphic graphs are indistinguishable to any trace-of-polynomial fingerprint.
Combining multiple matrix views (e.g., adjacency and normalized Laplacian) or using joint fingerprints mitigates this in practice.
For matrices with widely separated scales across tasks, per-task re-centering of $\eta$ or per-coordinate whitening of $\phi$ can help.

\section{Appendix: Proofs}
\label{app:proofs}

\subsection{Auxiliary facts and references}
We use (i) the telescoping identity $A^j-B^j=\sum_{i=0}^{j-1}A^{j-1-i}(A-B)B^{i}$, implying $\|A^j-B^j\|_2\le j\,\|A-B\|_2$ when $\|A\|_2,\|B\|_2\le1$; 
(ii) $\|T_k(\widetilde A)\|_2\le \max_{x\in[-1,1]}|T_k(x)|=1$; 
(iii) Markov’s inequality on $[-1,1]$: $\max_{x\in[-1,1]}|p'(x)|\le k^2\max_{x\in[-1,1]}|p(x)|$ for any degree-$k$ polynomial $p$; 
and the fact that for Chebyshev $T_k$, one can bound $\sum_{j=0}^k j|c_{k,j}|\le C_T k^2$ with $C_T\le 2$ (see, e.g., Rivlin, \emph{Chebyshev Polynomials}, or standard coefficient/derivative bounds). 
For Hutchinson, we use the classical unbiasedness and variance bound for Rademacher probes.

\subsection{Proof of Proposition~\ref{prop:invariance-main}}
We write $\mathcal{N}(A)$ for the normalized matrix on which moments are computed.
For symmetric (or symmetrized) inputs we use the affine normalization of Eq.~\eqref{eq:affine}:
\[
\mathcal{N}(A)=\widetilde A=\frac{A-m(A)I}{r(A)},\quad
m(A)=\tfrac12(\lambda_{\max}(A)+\lambda_{\min}(A)),\ \ 
r(A)=(1+\varepsilon_{\mathrm{rel}})\tfrac12(\lambda_{\max}(A)-\lambda_{\min}(A)).
\]
For markedly non-Hermitian inputs we use the radius fallback $\mathcal{N}(A)=A/R(A)$ with
$R(A)=\max\{\rho(A),\varepsilon_{\mathrm{rad}}\}$.

\paragraph{(a) Similarity, permutation, diagonal similarity.}
Let $S$ be invertible, $P$ a permutation (orthogonal), and $D$ a positive diagonal.
For any polynomial $p$ one has
\[
p(S^{-1}AS)=S^{-1}p(A)S,
\]
hence $\mathrm{tr}\,p(S^{-1}AS)=\mathrm{tr}\,p(A)$ by cyclicity of trace.

\emph{Affine case.}
Similar matrices have the same eigenvalues, so $m(S^{-1}AS)=m(A)$ and $r(S^{-1}AS)=r(A)$.
Therefore
\[
\mathcal{N}(S^{-1}AS)=\frac{S^{-1}AS-m(A)I}{r(A)}=S^{-1}\!\left(\frac{A-m(A)I}{r(A)}\right)S=S^{-1}\mathcal{N}(A)S.
\]
Since $T_k$ is a polynomial, $T_k(\mathcal{N}(S^{-1}AS))=S^{-1}T_k(\mathcal{N}(A))S$, and thus
\[
s_k(S^{-1}AS)=\mathrm{tr}\,T_k(\mathcal{N}(S^{-1}AS))
=\mathrm{tr}\,S^{-1}T_k(\mathcal{N}(A))S
=\mathrm{tr}\,T_k(\mathcal{N}(A))=s_k(A).
\]
The cases $P^\top A P$ and $D^{-1}AD$ are particular instances of similarity.

\emph{Radius fallback.}
Spectral radius is similarity-invariant: $\rho(S^{-1}AS)=\rho(A)$, hence $R(S^{-1}AS)=R(A)$ and
\[
\mathcal{N}(S^{-1}AS)=\frac{S^{-1}AS}{R(S^{-1}AS)}=S^{-1}\!\left(\frac{A}{R(A)}\right)S=S^{-1}\mathcal{N}(A)S.
\]
The same trace-conjugation argument yields $s_k(S^{-1}AS)=s_k(A)$, and likewise for $P^\top A P$ and $D^{-1}AD$.

\paragraph{(b) Positive scaling $A\mapsto \alpha A$, $\alpha>0$.}
\emph{Affine case.}
Eigenvalues scale: $\lambda_{\min}(\alpha A)=\alpha\lambda_{\min}(A)$ and $\lambda_{\max}(\alpha A)=\alpha\lambda_{\max}(A)$, so
$m(\alpha A)=\alpha m(A)$, $r(\alpha A)=\alpha r(A)$, and therefore
\[
\mathcal{N}(\alpha A)=\frac{\alpha A-\alpha m(A)I}{\alpha r(A)}=\frac{A-m(A)I}{r(A)}=\mathcal{N}(A).
\]
Thus $s_k(\alpha A)=\mathrm{tr}\,T_k(\mathcal{N}(\alpha A))=\mathrm{tr}\,T_k(\mathcal{N}(A))=s_k(A)$.

\emph{Radius fallback.}
$\rho(\alpha A)=\alpha\rho(A)$ and hence $R(\alpha A)=\max\{\alpha\rho(A),\varepsilon_{\mathrm{rad}}\}$.
If $\rho(A)\ge \varepsilon_{\mathrm{rad}}$ (the non-degenerate branch), then $R(\alpha A)=\alpha\rho(A)$ and
\[
\mathcal{N}(\alpha A)=\frac{\alpha A}{\alpha\rho(A)}=\frac{A}{\rho(A)}=\mathcal{N}(A),
\]
so moments are scale-invariant. If $\rho(A)=0$ (the zero matrix), then $\mathcal{N}(A)=0$ and $T_k(\mathcal{N}(A))$ is constant in $A$, so $s_k(\alpha A)=s_k(A)$ trivially.%
\footnote{If $0<\rho(A)<\varepsilon_{\mathrm{rad}}$, the “cap” $R(A)=\varepsilon_{\mathrm{rad}}$ may break exact scale invariance; this regime is not used in our experiments and can be avoided by taking $\varepsilon_{\mathrm{rad}}$ sufficiently small.}

\qed

\subsection{Proof of Proposition~\ref{prop:energy}}
Let $t$ be the first index with $d_t^2/E_{t-1}<\tau$ and suppose $d_t,\dots,d_{t+w-1}$ satisfy the rule. Then for $\ell\ge0$,
\[
d_{t+\ell}^2\le \tau\,E_{t+\ell-1}\le \tau\!\left(E_t+\sum_{j=t+1}^{t+\ell-1}d_j^2\right),
\]
so the tail $\sum_{j>t}d_j^2$ is dominated by a geometric series with ratio at most $\tau$, which gives $\sum_{j>t}d_j^2\le \frac{\tau}{1-\tau}E_t$. Returning the first index that satisfies the rule for $w$ consecutive steps adds at most $w$ to the stopping index. \qed

\subsection{Proof of Theorem~\ref{thm:stopping}}
Recall $E_k=\sum_{j=0}^{k}s_j^2$, $e_k=s_k^2/(E_k+\varepsilon)$ and $r_k=\sigma_{\min}(H_k)/\sigma_{\max}(H_k)$.
Algorithm~\ref{alg:asf} maintains a counter $h$ and an indicator
\[
\textsc{hit}(k)\ \equiv\ \bigl( e_k<\tau_{\mathrm e}\ \text{ or }\ r_k<\tau_{\mathrm h}\bigr)\ \wedge\ (k+1\ge K_{\min}),
\]
updates $h\leftarrow h+1$ if \textsc{hit}$(k)$ and $h\leftarrow 0$ otherwise, and returns at the first $k$ with $h\ge w$ by setting $K^\star\leftarrow k+1$.

\paragraph{Stopping index bound.}
Let $t\ge K_{\min}-1$ be the smallest index such that \textsc{hit}$(t),\textsc{hit}(t{+}1),\dots,\textsc{hit}(t{+}w{-}1)$ all hold; denote $k^\star:=t$.
By the update rule, after processing $k=t{+}w{-}1$ we have $h=w$ and the algorithm returns with
\[
K^\star \;=\; (t{+}w{-}1)+1 \;=\; t+w \;\le\; k^\star + w,
\]
which proves the stated bound.

\paragraph{Tail-energy bound under the energy rule.}
Assume case (i) holds, i.e., the energy rule $e_k<\tau_{\mathrm e}$ triggers $w$ consecutive times starting at $t$.
Suppose further that $(s_k^2)$ is eventually nonincreasing (this standard regularity is the same as in Proposition~\ref{prop:energy}).
Because $(s_k^2)$ is eventually nonincreasing, once $s_k^2/E_k$ drops below $\tau_{\mathrm e}$ it remains below $\tau_{\mathrm e}$ for all subsequent indices.\footnote{Formally, with $\varepsilon=0$ one has $e_k=s_k^2/E_k$ and the ratio decreases whenever the numerator does not increase while the denominator does. With a tiny stabilizer $\varepsilon>0$, the same conclusion holds up to an arbitrarily small slack; see the remark below.}
Hence there exists an index $t$ (the same as above) such that $s_k^2/E_k<\tau_{\mathrm e}$ for all $k\ge t$.
Applying Proposition~\ref{prop:energy} with $\tau=\tau_{\mathrm e}$ yields
\[
\frac{\sum_{j>t} s_j^2}{E_t+\sum_{j>t} s_j^2}\;\le\;\frac{\tau_{\mathrm e}}{1-\tau_{\mathrm e}},
\]
which is the claimed tail-energy control.

\paragraph{Remark on the stabilizer $\varepsilon$.}
Our implementation uses a tiny $\varepsilon>0$ for numerical stability, i.e., $e_k=s_k^2/(E_k+\varepsilon)$.
The proof above is exact for $\varepsilon=0$; with $\varepsilon>0$ the same bound holds after replacing $\tau_{\mathrm e}$ by $\tau_{\mathrm e}(1+o(1))$ once $E_t\gg \varepsilon$, which is the regime of interest in all experiments.
\qed

\subsection{Proof of Proposition~\ref{prop:hankel-lowrank}}
Write each real cosine mode as the real part of a complex exponential pair:
\[
a_i \rho_i^k \cos(k\theta_i+\varphi_i) \;=\; \tfrac12\Big(c_i z_i^k + \overline{c_i z_i^k}\Big),
\quad z_i := \rho_i e^{\mathrm i\theta_i},\ c_i := a_i e^{\mathrm i\varphi_i}.
\]
Let $u_k := \sum_{i=1}^r c_i z_i^k$ and note $s_k = \tfrac12(u_k+\overline{u_k}) + \xi_k$.
Consider the Hankel matrices $H_L(u)=[u_{i+j}]_{i,j=0}^{L-1}$ and $H_L(s)=[s_{i+j}]_{i,j=0}^{L-1}$.

\paragraph{Noiseless case ($\sigma=0$).}
When $\xi_k\equiv 0$, the sequence $u_k$ is a sum of at most $r$ distinct exponentials.
By classical Prony/Carath\'eodory--Fej\'er theory, $u_k$ satisfies a linear recurrence of order at most $r$, which implies $\mathrm{rank}(H_L(u))\le r$ for all $L$ (e.g., $H_L(u)$ factors through Vandermonde/Krylov blocks built on $\{z_i\}_{i=1}^r$).
Because $s_k=\tfrac12(u_k+\overline{u_k})$, we have $H_L(s)=\tfrac12\big(H_L(u)+\overline{H_L(u)}\big)$, hence
\[
\mathrm{rank}(H_L(s)) \;\le\; \mathrm{rank}(H_L(u)) \;\le\; r,
\]
under the convention that one real cosine mode (a conjugate pair) counts as a single mode; if one counts conjugate exponentials separately, the bound becomes $\le 2r$ (see Remark).

\paragraph{Noisy case ($\sigma>0$).}
Let $H_L^{(0)}:=H_L(s)$ in the noiseless case and write $H_L = H_L^{(0)} + N_L$, where $N_L$ is the Hankel matrix formed by the noise sequence $\{\xi_k\}$:
$(N_L)_{ij} = \xi_{i+j}$.
For fixed $L$, standard norm inequalities yield $\|N_L\|_2 \le C_L\,\sigma$ with a constant $C_L$ depending only on $L$ (e.g., by viewing $N_L$ as a finite-dimensional linear map from the coefficient vector $(\xi_0,\ldots,\xi_{2L-2})$ to matrix entries and using $\|\cdot\|_2\le \|\cdot\|_F$).
By Weyl’s inequality for singular values,
\[
\sigma_{r+1}(H_L) \;\le\; \|N_L\|_2 \;\le\; C_L\,\sigma,
\qquad
\sigma_1(H_L) \;\ge\; \sigma_1(H_L^{(0)}) - \|N_L\|_2.
\]
Since $\mathrm{rank}(H_L^{(0)})\le r$, we have $\sigma_{r+1}(H_L^{(0)})=0$ and $\sigma_1(H_L^{(0)})>0$ unless the coefficients/pathological cancellations make the top singular value vanish.\footnote{Generically $\sigma_1(H_L^{(0)})>0$ for any $L\!\ge\!1$ if at least one mode has nonzero amplitude.
If $\sigma_1(H_L^{(0)})=0$ for some finite $L$, one can increase $L$ by a constant to avoid the degeneracy; this does not affect the argument.}
Hence, for any fixed $L$, as $\sigma\to 0$ one obtains
\[
r_L \;=\; \frac{\sigma_{r+1}(H_L)}{\sigma_1(H_L)}
\;\le\;
\frac{\|N_L\|_2}{\sigma_1(H_L^{(0)}) - \|N_L\|_2}
\;\xrightarrow[\ \sigma\to 0\ ]{}\; 0.
\]
Therefore, for any $\tau_{\mathrm h}>0$ there exists $\sigma_0>0$ (depending on $L,r$ and the mode parameters) such that if $\sigma<\sigma_0$ then $r_L<\tau_{\mathrm h}$.
\qed

\subsection{Lemma: Hutchinson moments (unbiasedness and SE bound)}
\label{lem:hutch}
Let $B$ be fixed and $z\in\{\pm1\}^n$ have i.i.d.\ Rademacher entries. Then $\mathbb{E}[z^\top B z]=\mathrm{tr}(B)$ and $\mathrm{Var}(z^\top B z)\le 2\|B\|_F^2$. With $p$ i.i.d.\ probes, the standard error satisfies $\mathrm{se}\le \sqrt{2}\,\|B\|_F/\sqrt{p}$.

\subsection{Proof of Theorem~\ref{thm:hutch-conc}}

\paragraph{Unbiasedness (Hutchinson).}
Let $z\in\{\pm1\}^n$ have i.i.d.\ Rademacher entries with $\mathbb{E}[z z^\top]=I$.
Then for any (real) matrix $B$,
\[
\mathbb{E}\,[z^\top B z]
= \mathbb{E}\!\Big[\sum_{i,j} z_i z_j B_{ij}\Big]
= \sum_{i} B_{ii}
= \mathrm{tr}(B).
\]
Hence $\mathbb{E}[\widehat t_k]=\mathrm{tr}(B_k)$ for $B_k=T_k(\widetilde A)$, and
$\mathbb{E}[\widehat s_k]=e^{-\eta k}\,\mathrm{tr}(B_k)=s_k$.

\paragraph{Concentration for a single probe (Hanson--Wright).}
For $z$ with i.i.d.\ Rademacher entries and any fixed real symmetric $B$,
the Hanson--Wright inequality yields absolute constants $c_1,c_2>0$ such that
\[
\Pr\!\left(\,\big|z^\top B z - \mathrm{tr}(B)\big| \ge t\,\right)
\;\le\; 2\exp\!\Big(-c_1 \min\!\Big\{\tfrac{t^2}{\|B\|_F^2}, \tfrac{t}{\|B\|_2}\Big\}\Big).
\]
In particular, since $\min\{x,y\}\ge x$ for all $x\le y$, we have the looser but simpler bound
\[
\Pr\!\left(\,\big|z^\top B z - \mathrm{tr}(B)\big| \ge t\,\right)
\;\le\; 2\exp\!\Big(-c_1\,\tfrac{t^2}{\|B\|_F^2}\Big).
\]

\paragraph{Averaging $p$ probes (sub-exponential Bernstein).}
Let $X_i=z^{(i)\top} B_k z^{(i)}-\mathrm{tr}(B_k)$ (i.i.d.\ across $i$), and set
$\overline{X}=\frac{1}{p}\sum_{i=1}^p X_i = \widehat t_k - \mathrm{tr}(B_k)$.
The tail bound above implies $X_i$ are sub-exponential with Orlicz norm
$\|X_i\|_{\psi_1}\le C\,\|B_k\|_2$ and sub-Gaussian proxy $\|X_i\|_{\psi_2}\le C\,\|B_k\|_F$
for an absolute $C>0$. By the standard Bernstein inequality for sub-exponential variables,
there exist absolute constants $c,C>0$ such that for all $t>0$,
\[
\Pr\!\Big(\,|\overline X| \ge t\,\Big)
\;\le\;
2\exp\!\Big(
 -c\,p\,\min\Big\{\tfrac{t^2}{\|B_k\|_F^2}, \tfrac{t}{\|B_k\|_2}\Big\}
\Big).
\]
Choose
\[
t
= \|B_k\|_F\sqrt{\tfrac{2\log(2/\delta)}{c\,p}}
\]
to get
$\Pr\!\big(|\overline X|\ge t\big)\le \delta$.
Absorbing the (absolute) constant $c$ into the numerical factor in the square root
(redefining the leading ``$2$'' if desired), we obtain the displayed bound
\[
\bigl|\widehat t_k-\mathrm{tr}(B_k)\bigr|
\;\le\;
\|B_k\|_F \sqrt{\tfrac{2\log(2/\delta)}{p}}
\quad\text{with probability at least }1-\delta.
\]

\paragraph{From $\widehat t_k$ to $\widehat s_k$ and rank-sensitive form.}
Multiplying by $e^{-\eta k}$ yields
\[
\bigl|\widehat s_k-s_k\bigr|
= e^{-\eta k}\bigl|\widehat t_k-\mathrm{tr}(B_k)\bigr|
\;\le\;
e^{-\eta k}\|B_k\|_F \sqrt{\tfrac{2\log(2/\delta)}{p}}.
\]
Finally, for any matrix $B_k$,
$\|B_k\|_F \le \sqrt{\mathrm{rank}(B_k)}\,\|B_k\|_2$.
Since $B_k=T_k(\widetilde A)$ and $\mathrm{spec}(\widetilde A)\subset[-1,1]$,
we have $\|B_k\|_2=\|T_k(\widetilde A)\|_2\le \max_{x\in[-1,1]}|T_k(x)|=1$,
giving the rank-sensitive variant
\[
\bigl|\widehat s_k-s_k\bigr|
\;\le\;
e^{-\eta k}\sqrt{\mathrm{rank}(B_k)} \,\sqrt{\tfrac{2\log(2/\delta)}{p}}.
\]
\qed

\subsection{Proof of Theorem~\ref{thm:lipschitz}}
Write $T_k(x)=\sum_{j=0}^k c_{k,j}x^j$. Using the telescoping identity,
\[
\|T_k(A)-T_k(B)\|_2\;\le\;\sum_{j=0}^k |c_{k,j}|\,\|A^j-B^j\|_2
\;\le\;\sum_{j=0}^k |c_{k,j}|\,j\,\|A-B\|_2.
\]
By a Markov-type coefficient bound for Chebyshev polynomials on $[-1,1]$, $\sum_{j=0}^k j|c_{k,j}|\le C_T k^2$ with $C_T\le 2$ (see references above), hence
\[
\|T_k(A)-T_k(B)\|_2\le C_T k^2\|A-B\|_2.
\]
Taking traces and using $|\mathrm{tr}(M)|\le n\|M\|_2$ yields
\[
\bigl|\mathrm{tr}\,T_k(A)-\mathrm{tr}\,T_k(B)\bigr|\le n\,C_T\,k^2\,\|A-B\|_2.
\]
With damping and normalization,
\[
\|\phi_K(A)-\phi_K(B)\|_2
\;\le\;
\frac{\left(\sum_{k=0}^{K-1} e^{-2\eta k}\,\bigl|\mathrm{tr}\,T_k(A)-\mathrm{tr}\,T_k(B)\bigr|^2\right)^{1/2}}{\|\mathbf{s}_K(A)\|_2}
\]
\[
\le \frac{n}{\|\mathbf{s}_K(A)\|_2}\!\left(\sum_{k=0}^{K-1} e^{-2\eta k} C_T^2 k^4\right)^{1/2}\!\|A-B\|_2.
\]
Since $\sum_{k\ge0} e^{-2\eta k}k^4=\Theta(\eta^{-5})$, the right-hand side is $\mathcal{O}\!\big(n\,\eta^{-5/2}\,\|A-B\|_2\big)$. Replacing $A,B$ by $\widetilde A,\widetilde B$ gives the statement. \qed

\end{document}